\documentclass[12pt]{amsart}
\usepackage{amsmath,amstext,amsfonts,
amssymb}
\usepackage{eucal}
\usepackage{diagrams}                                           %
\diagramstyle[scriptlabels,PostScript=dvips]               %
\newarrow{Dash}{}{dash}{}{dash}{}                          %
\renewcommand{\subsubsection}[1]{\addtocounter{subsubsection}{1}
{\ \\[3pt]\bf \thesubsubsection. \  #1} }

\swapnumbers

\newtheorem{lem}[subsubsection]{Lemma}
\newtheorem{prp}[subsubsection]{Proposition}
\newtheorem{crl}[subsubsection]{Corollary}

\newtheorem{Prp}[subsection]{Proposition}

\newtheorem*{thm*}{Theorem}
\newtheorem*{lem*}{Lemma}
\newtheorem*{prp*}{Proposition}
\newtheorem*{crl*}{Corollary}

{  \theoremstyle{definition}
           \newtheorem{dfn}[subsubsection]{Definition}

           \newtheorem{Dfn}[subsection]{Definition}

           \newtheorem*{dfn*}{Definition}
           \newtheorem*{exm*}{Example}
           \newtheorem*{rem*}{Remark}
           \newtheorem*{rems*}{Remarks}
}

\newcommand{\holim}{\operatorname{holim}}

\newcommand{\p}{\partial}
\newcommand{\wt}{\widetilde}

\newcommand{\Mor}{\mathrm{Mor}}

\newcommand{\cC}{\mathcal{C}}
\newcommand{\cD}{\mathcal{D}}

\newcommand{\cH}{\mathcal{H}}

\newcommand{\cN}{\mathcal{N}}
\newcommand{\cO}{\mathcal{O}}

\newcommand{\cW}{\mathcal{W}}

\newcommand{\id}{\mathrm{id}}
\newcommand{\pr}{\mathrm{pr}}

\newcommand{\N}{\mathbb{N}}

\newcommand{\bW}{\bar{W}}
\newcommand{\bbW}{\mathbb{W}}
\newcommand{\bE}{\mathbb{E}}
\newcommand{\wh}{\widehat}

\newcommand{\Hom}{\operatorname{Hom}}

\newcommand{\Def}{\mathrm{Def}}
\newcommand{\sDef}{\cD\mathit{ef}}

\newcommand{\Ind}{\mathrm{Ind}}

\newcommand{\fN}{\mathfrak{N}}
\newcommand{\fl}{{\mathit{f}\ell}}

\newcommand{\LIM}{\mathrm{LIM}}

\newcommand{\op}{\mathrm{op}}

\newcommand{\Set}{\mathtt{Set}}
\newcommand{\sSet}{\mathtt{sSet}}
\newcommand{\simpl}{\mathtt{sSet}}
\newcommand{\sCat}{\mathtt{sCat}}
\newcommand{\Cat}{\mathtt{Cat}}
\newcommand{\CAT}{\mathbf{\Delta}_\fN}
\newcommand{\SC}{\mathbf{\Delta}_\cN} 
\newcommand{\WC}{\mathbf{\Delta}_{\bar{W}}}

\newcommand{\shom}{\cH\!\mathit{om}}

\newcommand{\Ob}{\operatorname{Ob}}

\begin{document}
\title[]{Homotopy coherent nerve in deformation theory}
\author{Vladimir Hinich}    
\address{Department of Mathematics, University of Haifa,    
Mount Carmel, Haifa 31905,  Israel} 
\email{hinich@math.haifa.ac.il}    
\begin{abstract}
The main object of the note is to fix an earlier error
of the author, \cite{dgc}, claiming that 
the (standard) simplicial nerve preserves fibrations of simplicially 
enriched categories.
This erroneous claim was used by the author in his  deformation
theory constructions in~\cite{dgc,dha,dsa}. 

Fortunately, the problem disappears if one replaces the standard
simplicial nerve with another one, called {\em homotopy coherent
nerve}. 

In this note we recall the definition of homotopy coherent nerve and 
prove some its properties necessary to justify the papers~\cite{dgc,dha,dsa}.
We present as well a generalization of the notion of fibered
categories which is convenient once one works with enriched categories.

\end{abstract}
\maketitle    

\setcounter{section}{-1}
\section{Introduction}

In a series of papers~\cite{dgc,dha,dsa} dealing with formal deformation 
theory the author used
a non-standard model category structure for the category $\sCat$ of 
simplicial (more precisely, simplicially enriched) categories, and a 
standard simplicial nerve functor
$$\cN:\sCat\rTo\simpl.$$
An important claim, Proposition A3.2 of~\cite{dgc}, (which was found to be
wrong by Sergey Mozgovoy) claimed that the nerve functor preserves
fibrations. Preservation of fibrations is crucial in \cite{dgc,dha,dsa}
at two points. The less important among the two is connected to a simplicial
structure on $\sCat$ which we believed (based on the preservation of 
fibrations) to provide a simplicial framing in the sense of~\cite{hir}, 16.6,
see Section~\ref{ss:holim}.
A more important gap appears in the comparison of two versions of a 
deformation functor: one with values in simplicial categories, and the other
(obtained by application of the nerve functor to the first one) with values
in simplicial sets. We used preservation of fibrations for the nerve functor
to deduce that both versions of the deformation functor satisfy a descent 
property.

Fortunately, it turns out that another version of nerve functor,
$$\fN:\sCat\rTo\simpl,$$
known to homotopy theorists as {\em homotopy coherent nerve},
satisfies the required properties.


When applied to a fibrant groupoid (a fibrant simplicial category $\cC$
whose $0$-arrows are invertible up to homotopy), the homotopy coherent 
nerve $\fN$ is homotopy equivalent to the standard nerve $\cN$. 
As well, $\fN$ preserves fibrations of fibrant groupoids. In this way most 
of the claims made in~\cite{dgc,dha,dsa} become true with $\cN$ replaced 
with the new nerve.

The homotopy coherent nerve functor does not help in providing a simplicial
framing for $\sCat$; thus, we have to use a choice of simplicial framing
to deal with homotopy limits which abound in~\cite{dgc,dsa}. To cope with
weak functors, we replace the pseudofunctors defined as in \cite{dsa}, A1.11,
with a more convenient notion of quasifibered categories, see~\ref{ss:QF}. 
The notion of quasifibered category may be also convenient in other contexts.
{It is highly probable that this notion is well-known.}  
One of the obvious conveniences of this notion is
a possibility to apply Dwyer-Kan localization ``fiberwise'' to a quasifibered 
category.

The notion of homotopy coherent nerve belongs apparently to Cordier
\cite{cordier1}. The adjoint functor $\CAT$ applied to
the nerve $N\cC$ of a category $\cC$, provides a ``simplicial thickening'' 
of $\cC$ used by Cordier to define homotopy coherent diagrams in the simplicial
context. Theorem 2.1 of~\cite{corpor} claims that the homotopy coherent
nerve of a fibrant simplicial category is weakly Kan, that is satisfies 
the right lifting property with respect to all {\em inner horns}
$\Lambda^n_i\rTo\Delta^n,\ i=1,\ldots,n-1.$
The latter result is a small part of J.~Lurie's theorem
1.3.4.1 from~\cite{lurie}, 
claiming that the pair of functors $\fN$ and $\CAT$ define a Quillen 
equivalence between the homotopy categories of $\sCat$ and of the category of
simplicial sets endowed with the Joyal model structure. The same result is
also announced by Joyal~\cite{joyal}.

In hindsight, one can explain the relevance of homotopy coherent nerve
for deformation theory as follows. In the version of deformation theory
studied in~\cite{dgc,dha,dsa} the deformation functors assign to an artinian
commutative ring an $\infty$-groupoid.
Kan simplicial sets and  fibrant groupoids are two
ways to formalize this notion. In the same spirit weakly Kan simplicial
sets and fibrant simplicial categories both formalize the notion of
$\infty$-category in terms of Joyal and Lurie. The equivalence between 
the two is established by the adjoint pair $(\CAT,\fN)$. Thus, it is very 
natural to use the same functors to connect between $\sCat$- and 
$\simpl$-valued versions of deformation functors.

\subsection{Overview}

The standard nerve functor
$$\cN:\sCat\rTo\simpl$$
is represented by a cosimplicial object $\SC^\bullet$ in $\sCat$
such that the objects of $\SC^n$ are the natural numbers $0,\ldots,n$
and the morphisms of $\SC^n$ are freely generated by
the simplicial sets $\Hom(i-1,i)$ which are standard $n$-simplices.

The homotopy coherent nerve is represented by another cosimplicial object 
which we denote
$\CAT^\bullet$ (see definition in~\ref{ss:CAT-1}, \ref{ss:CAT-2}). 
Without giving here a precise definition, let us just indicate that for
a simplicial category $\cC$ the $0$-simplices of the new nerve $\fN(\cC)$
are just the objects of $\cC$, the $1$-simplices are $0$-arrows in $\cC$, and
two-simplices are the collections of three objects $0,1,2$, three
$0$-arrows $f_{ij}:i\to j,\ 0\leq i<j\leq 2$, and a $1$-arrow from 
$f_{12}\circ f_{01}$ to $f_{02}$.

The homotopy coherent nerve functor $\fN:\sCat\rTo\simpl$ enjoys the 
following properties.

\begin{prp*}
\begin{itemize}
\item[(1)](Corollary~\ref{N-vs-N}) There is a natural map 
$\cN(\cC)\rTo\fN(\cC)$ which is homotopy equivalence for fibrant groupoids.
\item[(2)](Proposition~\ref{prp:wef})
The homotopy coherent nerve functor $\fN$ preserves weak equivalences of 
fibrant simplicial categories.
\item[(3)](Proposition~\ref{prp:fibrations}, see also~\ref{proof:fibrations})
The functor $\fN$ preserves fibrations of fibrant groupoids.
\end{itemize}
\end{prp*}

Our central result is the third claim of the proposition proven in
~\ref{proof:fibrations}. The claims 1 and 2 of the proposition can be probably
deduced from \cite{corpor}, \cite{lurie}, 1.3.4.1, and~\cite{joyal}.
We present their proof for convenience.


The above mentioned properties of $\fN$ allows one to deduce that $\fN$ 
commutes with homotopy limits of (weak) functors with 
values in fibrant groupoids --- see precise formulation in~\ref{N-holim}.

\subsection{Plan} In Section~\ref{sec:models} we recall two model category 
structures on $\sCat$, the standard one due to W.~Dwyer, P.~Hirschhorn. D.~Kan
 and J.~Bergner \cite{Bergner}, 
and the one we defined in~\cite{dgc}. The standard model
structure is more natural; we leave the non-standard one mostly for the sake
of notational compatibility. 
\footnote{We still feel the non-standard model structure
may be of some use.}
Throughout the paper, we use both of them. In Section~\ref{sec:new-nerve} 
we present a definition of  the  homotopy coherent
nerve and prove some properties (more technical proofs are moved to
Section~\ref{sec:details}). We define quasifibered categories and prove some 
their properties in Section~\ref{sec:qf}. 
We present in~\ref{ss:examples} two counterexamples
dealing with preservation of fibrations by the nerve functors.

Finally, in Section~\ref{sec:corr} we present the list of corrections for
\cite{dgc,dha,dsa}.

{\bf Acknowledgement.} \ The author is deeply grateful to Sergey Mozgovoy
who has pointed out to an error in~\cite{dgc}.
Without him this note could not have been written.  The author is 
also grateful to the referee whose remarks helped to improve the exposition.

\section{Model category structures}
\label{sec:models}

There are at least two model category structures on $\sCat$. The first one,
the more standard and more natural among the two,
was suggested by Dwyer, Hirschhorn and Kan in a preliminary version 
\footnote{now unavailable} of their book and described in a detail by 
J.~Bergner, \cite{Bergner}. Another structure defined in~\cite{dgc} has more
weak equivalences (someone would say, too weak equivalences), and the same 
cofibrations. Both structures are cofibrantly generated.

Let us recall the definitions.

\subsection{}
The functor
$$\pi_0:\sCat\to\Cat$$
is defined by the formulas 
\begin{eqnarray}
  \label{eq:pi_01}
  \Ob \pi_0(X)=\Ob X\\
  \label{eq:pi_02}
  \Hom_{\pi_0(X)}(x,y)=\pi_0(\Hom_X(x,y)).
\end{eqnarray}

\begin{Dfn}\label{scat-we}
A map $f:\cC\to\cD$ in $\sCat$ is called  {\em a weak equivalence}
if the map $\pi_0(f):\pi_0(\cC)\to\pi_0(\cD)$ induces a homotopy
equivalence of the nerves, and for each $x,y\in\Ob(\cC)$ the map
$$ \shom(x,y)\to\shom(f(x),f(y))$$
of the simplicial Hom-sets is a weak equivalence. 
\end{Dfn}

\begin{Dfn}\label{scat-e}
A weak equivalence $f:\cC\to\cD$ in $\sCat$ is called  {\em a strong 
equivalence}
if the map $\pi_0(f):\pi_0(\cC)\to\pi_0(\cD)$ is an equivalence of
categories.

\end{Dfn}

\subsection{}
Cofibrations in $\sCat$ are generated by the following maps:

\begin{itemize}
\item[(i)] $\emptyset\to *$, the functor from the empty simplicial 
category to a one-point category.
\item[(ii)] For each  cofibration $K\to L$ in $\simpl$ the induced  map from
$K_{01}$ to $L_{01}$. Here $K_{01}$ denotes the simplicial category having 
two objects $0$ and $1$ with the only nontrivial maps $K=\shom(0,1)$.
\end{itemize}

The notions of cofibration and of (weak) equivalence define two model category
structures on $\sCat$.

\begin{Prp}
Cofibrations and strong equivalences as above define a model structure on 
$\sCat$ (the standard model structure).
\end{Prp}
The fibrations in the standard model category structure will be called 
{\em weak fibrations}. They are characterized by the following proposition.

\begin{Prp}
\label{prop-weak-fibration}
A map $f:\cC\rTo\cD$ is a weak fibration if and only if the following
properties are satisfied.
\begin{itemize}
\item For all $x,y\in\cC$ the map $\shom_\cC(x,y)\to\shom_\cD(fx,fy)$
is a Kan fibration.
\item For all $x\in\cC$, any {\em homotopy equivalence} $\alpha:fx\to z$
in $\cD$ can be lifted to a homotopy equivalence $x\to y$ in $\cC$.
\end{itemize}
\end{Prp}

Here a $0$-map $\alpha$ is called homotopy equivalence if $\pi_0(\alpha)$
is an isomorphism.

\begin{Prp}
Cofibrations and weak equivalences as above define a model structure on 
$\sCat$ (the non-standard model structure).
\end{Prp}
The corresponding fibrations are called {\em strong fibrations}.

\subsection{}
The collection of weak acyclic cofibrations (that is, cofibrations which are
simultaneously weak equivalences) is generated by the following collection.

\begin{itemize}
\item[(i)] The maps $\partial^{0,1}:[0]\to[1]$
from the terminal category to the one-arrow category.  
\item[(ii)] For each  acyclic cofibration $K\to L$ in $\simpl$, the induced  
map from $K_{01}$ to $L_{01}$. 
Here, as above, $K_{01}$ is the simplicial category 
having two objects $0$ and $1$ with the only nontrivial maps $K=\shom(0,1)$.
\end{itemize}

A generating collection of strong acyclic cofibrations is given 
in~\cite{Bergner}. 

\subsection{}
A simplicial category $\cW$ is called {\em weak groupoid} if $\pi_0(\cW)$
is a groupoid. Equivalently, weak groupoids are simplicial categories
strongly equivalent to simplicial groupoids.

Note that for a simplicial category $\cC$ the notions of weakly fibrant
and strongly fibrant category coincide. We will call them simply
{\em fibrant categories}.

Weak groupoids are not necessarily fibrant. A fibrant weak groupoid
will be called simply {\em fibrant groupoid}. Of course, simplicial groupoids
are fibrant.

\section{Homotopy coherent nerve}
\label{sec:new-nerve}

For each $n\in\N$ we will define a simplicial category $\CAT^n$.
The collection $\{n\mapsto\CAT^n\}$ gives a cosimplicial object in $\sCat$,
and the nerve functor
$$\fN:\sCat\rTo\sSet$$
will be defined by the formula
$$\fN(\cC)_n=\Hom(\CAT^n,\cC).$$

\subsection{Notation}
We denote by $\N$ the set of nonnegative integers. 
The objects of the standard category of simplices $\Delta$ are denoted 
$[n]$ or just $n$ according to convenience. 
For $a,b\in\N$ we denote by $(a,b)$ the segment $\{c\in\N|a<c<b\}$.
This notation is only used when $a\leq b$.

Some notation for categories. $\Set$ is the category of sets, $\simpl$ 
that of simplicial sets and $\sCat$ that of simplicial (=simplicially enriched)
categories.

\subsection{The simplicial categories $\CAT^n$}
\label{ss:CAT-1}

The simplicial category $\CAT^n$ playing the role of a standard $n$-simplex
in $\sCat$ is a cofibrant resolution of the standard $n$-simplex in $\Cat$
defined by the ordered set $[n]=\{0,\ldots,n\}$.

Thus, $\Ob\CAT^n=\{0,\ldots,n\}$. For $a\leq b$ in $\Ob\CAT^n$ 
we define $\CAT^n(a,b):=I^{(a,b)}$, 
the power of the standard $1$-simplex $I=\Delta^1$.

The composition 
$$\CAT^n(b,c)\times\CAT^n(a,b)\rTo\CAT^n(a,c)$$
is given by the map
$$ I^{(b,c)}\times I^{(a,b)}\rTo I^{(a,c)}$$
given by $\id^{(b,c)}\times\p^1\times\id^{(a,b)}$
where
$$ \p^i:[k]\rTo\left[k+1\right]$$
is the standard face map missing the value $i$
\footnote{the composition  map inserts $0$ at the place $b$.}.

The fact that $\pi_0(\CAT^n)=[n]$ is obvious; the simplicial category
$\CAT^n$ is cofibrant. The set of free generators for morphisms of
$\CAT^n$ is presented in~\ref{cat-n-are-cofibrant}.

\subsection{Cosimplicial structure on $n\mapsto\CAT^n$}
\label{ss:CAT-2}

Let $f:[m]\rTo\left[n\right]$ be a non-decreasing map. We will now
define a simplicial functor $f_*:\CAT^m\rTo\CAT^n$.

This will define a cosimplicial object $\CAT^\bullet$ in $\sCat$.
The functor $f_*$ is given by $f$ on the objects.
To specify the map $f_*:I^{(a,b)}\rTo I^{(f(a),f(b))}$ of
simplicial sets, we notice that the simplicial sets $I^X$, $X\in\Set$,
are nerves of posets; a map of such simplicial sets is uniquely determined
by its restriction on $0$-simplices.

 There exists a unique indecomposable $0$-simplex
in $\shom(a,b)$; we denote it by $\phi(a,b)$ and it is given
by the formula $\phi(a,b)_x=1$ for $x\in(a,b)$.
We write $\phi(a,a)=\id_a$ for simplicity. 

The map $f_*$ is defined by the formula 
$$ f_*(\phi(a,b))=\phi(f(a),f(b)).$$

\subsubsection{}
\label{functor-cat}
Define the functor 
$$ \CAT:\simpl\rTo\sCat$$
by the properties
\begin{itemize}
\item $\CAT$ commutes with colimits.
\item $\CAT(\Delta^n)=\CAT^n$.
\end{itemize}

\subsection{Properties of $\CAT$}
The following is a direct consequence of the definitions.

\begin{lem}
The map $\CAT(\p\Delta^n)\rTo\CAT^n$ is a cofibration.
\end{lem}
\qed

The functor $\CAT$ does not necessarily carry acyclic cofibrations
to weak acyclic cofibrations. A typical example is the map
$\CAT(\Lambda^n_0)\rTo\CAT(\Delta^n)$ induced by the zeroth horn
$\Lambda^n_0\rTo\Delta^n$. A calculation (see~\ref{lambda0n}) shows that 
the map
$$ \CAT(\Lambda^n_0)(1,n)\rTo\CAT(\Delta^n)(1,n)=I^{(1,n)}$$
identifies $ \CAT(\Lambda^n_0)(1,n)$ with the boundary of the cube
$I^{(1,n)}$.\footnote
{An even more obvious example $\CAT^0\rTo\CAT^1$ shows that $\CAT$ does 
not carry acyclic cofibrations to {\em strong} acyclic cofibrations.}

The following, however, takes place.
\begin{prp}
\label{prp:inner-horns}
For $i\ne 0,n$ the map
$$ \CAT(\Lambda^n_i)\rTo\CAT^n$$
is a strong acyclic cofibration.
\end{prp}

The proof is given in~\ref{inner-horns}.

\subsection{Homotopy coherent nerve and its properties}
\label{ss:newnerveprop}

We define the functor
$$\fN:\sCat\rTo\simpl$$
by the formula
$\fN_n(X)=\Hom_{\sCat}(\CAT^n,X).$

The functor $\fN$ is right adjoint to the functor
$$\CAT:\simpl\rTo\sCat$$
defined in~\ref{functor-cat}.

The functor $\fN$ preserves acyclic fibrations since $\CAT$ preserves 
cofibrations. It preserves cofibrations simply because 
injective maps are cofibrations in $\simpl$. 
One can not expect $\fN$ to preserve equivalences
and fibrations in general. We will prove this property of $\fN$ for 
special classes of simplicial categories.

\begin{prp}
\label{prp:wef}
$\fN$ preserves weak equivalences
of fibrant simplicial categories.
\end{prp}
\begin{proof}
Let $f:X\rTo Y$ be a weak equivalence with $X$ and $Y$ fibrant.
We have to check that for each $n\geq 0$ the induced map
$$ \pi_n(\fN(X))\rTo\pi_n(\fN(Y))$$
is an isomorphism.
We already know this if $f$ is an acyclic fibration. Thus, we can assume
that $f$ is a (weak) acyclic cofibration. Since $X$ is fibrant, $f$ splits: 
there exists $g:Y\rTo X$
such that $gf=\id_X$. Decompose $g=pi$ where $p:Z\rTo X$ is an acyclic 
fibration and $i:Y\rTo Z$ is a weak acyclic cofibration. Since $p$ is
an acyclic fibration, $Z$ is fibrant and $\pi_n(\fN(p))$ is an isomorphism.
 Thus, 
$$\pi_n(\fN(i))\circ\pi_n(\fN(f))$$
is a bijection which implies that $\pi_n(\fN(f))$ is injective. Since this 
is true for any weak acyclic cofibration $f$, this is also true for $i$.

Therefore $\pi_n(\fN(f))$ is necessarily bijective.
\end{proof}

\begin{prp}
\label{prp:fibrations}
The functor $\fN$ carries weak fibrations
\footnote{for fibrant groupoids they are automatically strong fibrations}
 of fibrant groupoids into Kan 
fibrations.
\end{prp}

The proof is given in~\ref{proof:fibrations}. Some counterexamples
to the preservation of fibrations by the nerve functor are given
in~\ref{ss:examples}.

\subsection{Comparison to other versions of nerve functor.}

In the appendix to~\cite{dgc} we used the standard definition of 
simplicial nerve
$$\cN:\sCat\rTo\sSet$$
as the diagonal of the bisimplicial set obtained by application the usual
nerve functor to the categories forming a simplicial category. The functor
$\cN$, similarly to $\fN$, is represented by a cosimplicial object 
$\SC^\bullet$ in $\sCat$ playing a role of the collection of standard 
simplices. 

Recall that $\Ob\SC^n=\{0,\ldots,n\}$ and that the morphisms of $\SC^n$
are freely generated by $n$-simplices $f_i\in\shom(i-1,i)$ for $i=1,
\ldots,n$, see~\cite{dgc}, A3.1.

Another variation of simplicial nerve functor we want to mention
is the functor $\bW$, see~\cite{may}, chapter 4. The functor $\bW$ is also 
represented by a cosimplicial object in $\sCat$ which we denote $\WC^\bullet$.
The simplicial category $\WC^n$
is defined as
\begin{itemize}
\item $\Ob\WC^n=\{0,\ldots,n\}$.
\item Norphisms of $\WC^n$ are freely generated by $n-i$-simplices
$g_i\in\shom(i-1,i)$, for $i=1,\ldots,n$.
\end{itemize}
The cosimplicial structure on the collection
$n\mapsto \WC^n$ is given by the formulas
\begin{itemize}
\item For $i\ne n-1$ \ $\p^i:\WC^{n-1}\rTo\WC^n$ carries
$g_j$ to $d_{i-j}g_j$ for $j<i$, to $g_{i+1}d_0g_i$ for $j=i$ and to
$g_{j+1}$ for $j>i$. 
\item $\p^{n-1}(g_j)=d_{n-1-j}g_j$.
\item $\sigma^i:\WC^n\rTo\WC^{n-1}$ carries $g_j$ to $s_{i-j}g_j$
for $j\leq i$, to $\id_j$ for $j=i+1$, and to $g_{j-1}$ for $j>i+1$.
\end{itemize}

In this subsection we show that the functors $\cN$ and $\fN$ agree up to 
homotopy on fibrant groupoids.

\subsubsection{Morphisms  $\CAT^\bullet\rTo\WC^\bullet\rTo\SC^\bullet$  }

The morphism $\pi:\WC^\bullet\rTo\SC^\bullet$ of cosimplicial objects in 
$\sCat$ inducing a homotopy equivalence of the nerves
$$ \cN(\cC)\rTo\bW(\cC)$$
for any simplicial groupoid $\cC$, is presented in~\cite{dgc}, A5.1.
It is identical on the objects, and sends $g_i$ to $d_0^if_i$ in the standard 
notation. 

We claim there exists a unique map $\tau:\CAT^\bullet\rTo\WC^\bullet$
bijective on the objects.
The map $\rho^1:\CAT^1\rTo\bW^1$ sends the unique morphism from $0$ to $1$
in $\CAT^1$ to the unique arrow in $\WC^1$. Let $0\leq a<b\leq n$.
The space $\CAT^n(a,b)$ contains a unique indecomposable
$0$-simplex denoted $\phi(a,b)$. Its image $\psi(a,b):=\tau(\phi(a,b))$
can be easily calculated. In fact, the map $d:\Delta^1\rTo\Delta^n$
defined by the formulas $d(0)=a,\ d(1)=b$, induces the maps
$d_*:\CAT^1\rTo\CAT^n$ and $d_*:\WC^1\rTo\WC^n$; since 
$\phi(a,b)=d_*\phi(0,1)$ one should have $\psi(a,b)=d_*\psi(0,1)$.
A simple calculation gives
\begin{equation}\label{psi(a,b)}
\psi(a,b)=\iota(g_b)\cdot \iota d_0(g_{b-1})\cdots \iota d_0^{b-a-1}(g_{a+1}),
\end{equation}
where $\iota=d_1^{n-b}$ assigns to each simplex its initial vertex.
Since any $0$-arrow in $\CAT^n$ is uniquely presented as a composition
of $\phi(a,b)$, the map $\rho$ is uniquely defined on $0$-arrows.

Note that the simplicial sets of morphisms both in $\CAT^n$ and in $\WC^n$
are nerves of posets; a map of such simplicial sets is uniquely defined
by its restriction on $0$-simplices which should be monotone.

Thus, it remains to check the monotonicity of the map defined
by the formula~(\ref{psi(a,b)}): one has to check that
$$\psi(b,c)\psi(a,b)<\psi(a,c).$$
This immediately follows from the inequality $d_1x<d_0x$ for any one-simplex 
$x$.

\begin{prp}\label{fn-vs-wb}Let $\cC$ be a simplicial groupoid. Then
the map $\bW(\cC)\rTo\fN(\cC)$ induced by $\tau$ is an equivalence.
\end{prp}
\begin{proof}
We can assume that $\cC$ is connected. Moreover, the claim easily reduces
to the case $\cC$ is a simplicial group (that is has only one object $*$).
It is known that $\pi_n(\bW(\cC))=\pi_{n-1}(\shom_\cC(*,*))$.
Let us calculate $\pi_n(\fN(\cC))$. By \ref{prp:fibrations} $\fN(\cC)$
is a Kan simplicial set. Thus, $\pi_n(\fN(\cC))$ can be easily calculated.
Let 
$$S^n=\Delta^n/\p\Delta^n,\quad 
D^{n+1}=\Delta^{n+1}/\cup_{i>0}\p^i(\Delta^n).$$
Then the face map $\p^0$ induces an embedding $S^n\rTo D^{n+1}$
and $\pi_n(\fN(\cC))$ is the quotient of the abelian group
$\Hom(S^n,\fN(\cC))=\Hom(\CAT(S^n),\cC)$ by the subgroup 
$\Hom(\CAT(D^{n+1}),\cC)$. The simplicial categories $\CAT(S^n)$ and 
$\CAT(D^{n+1})$ are explicitly described in~\ref{cat-examples}. 
A map $\CAT(S^n)\rTo\cC$ is given by a map
$$ \underbrace{S^1\wedge\ldots\wedge S^1}_{n-1} \rTo\shom_{\cC}(*,*)$$
of simplicial sets, whereas a map $\CAT(D^{n+1})\rTo\cC$ corresponds to
a map
$$ I\wedge \underbrace{S^1\wedge\ldots\wedge S^1}_{n-1} \rTo\shom_{\cC}(*,*).$$
Since $\shom_{\cC}(*,*)$ is a simplicial group, it is fibrant, and
the quotient gives the required result.
\end{proof}

\begin{crl}
\label{N-vs-N}
The composition map $\CAT^\bullet\rTo\SC^\bullet$ induces
a homotopy equivalence
\begin{equation}\label{nerve-eq}
 \cN(\cC)\rTo\fN(\cC)
\end{equation}
for any fibrant groupoid $\cC$. 
\end{crl}
\begin{proof}
Both $\fN$ and $\cN$ preserve weak equivalences of fibrant simplicial 
categories. Thus, the map~(\ref{nerve-eq}) is equivalence since it is
equivalence with $\cC$ replaced by its full hammock localization
 $L^H(\cC,\cC)$.
\end{proof}

\section{quasifibered categories and homotopy limits}
\label{sec:qf}

The simplicial structure on $\sCat$ defined in~\cite{dgc} does not
\footnote{contrary to the claim of~\cite{dgc}, A.4.3} provide a simplicial 
framing even for simplicial groupoids.

Thus, we need a choice of simplicial framing to construct homotopy limits.

In this section we choose a slightly different approach to 
dealing with weak functors with values in categories or simplicial categories.
In~\cite{dsa} we used the notion of a pseudofuctor which we now find 
inconvenient for two reasons. The first is the lack of simplicial model 
category structure on $\sCat$. The second inconvenience is connected to the
fact that, even though one can describe fibered categories in the language 
of pseudofunctors, this requires unnecessary choices to be made.


Let us, however, start with stict functors and simplicial framings.

\subsection{Homotopy limits}
\label{ss:holim}
To deal with homotopy limits of simplicial categories, we have to choose 
once and forever a simplicial framing in the standard model category 
structure. 
This means that a functorial map
$\cC\rTo\wt{\cC}_\bullet$ of simplicial objects in $\sCat$ is given
($\cC$ being considered as a constant simplicial object)
such that
\begin{itemize}
\item $\cC\rTo\wt{\cC}_0$ is an isomorphism.
\item $\cC\rTo\wt{\cC}_n$ is a strong equivalence for $\cC$ fibrant.
\item Define for $S\in\simpl$ a function space $\cC^S=
\lim_{[n]\to S}\wt{\cC}_n$. Then for fibrant $\cC$ any cofibration $S\to T$
defines a weak fibration $\cC^T\rTo\cC^S$.
\end{itemize}

Once we have a simplical framing,
we can define the homotopy limit of a functor $F:B^\op\rTo\sCat$
by the formula
\begin{equation}
\label{holimscat}
\holim(F)=\lim_{p\to q}(\wt{F}^q)^{\Delta^p},
\end{equation}
where the cosimplicial replacement $\wt{F}:\Delta\rTo\sCat$ is defined, as
usual, by the formula
$$ \wt{F}^n=\prod_{x_0\to\ldots\to x_n\in B}F(x_0).$$

Of course, the notion of a strict functor with values in $\sCat$ or even in
$\Cat$ is not very practical. We present below a convenient weak substitute.

\subsection{Quasifibered categories}
\label{ss:QF}
\subsubsection{}
Recall a few standard notions from \cite{sga1}, \'exp.~VI.

Let $p:E\rTo B$ be a category over $B$. An arrow $f:x\to y$ in $E$
over $\phi:b\to c$ in $B$ is called {\em cartesian} if the natural
map of Hom-sets
$$ \Hom_b(z,x)\rTo\Hom_\phi(z,y)$$
is a bijection. Here and below $\Hom_\phi$ denotes the set of arrows over 
$\phi$ and $\Hom_b$ is the same as $\Hom_{\id_b}$.

The category $p:E\rTo B$ over $B$ is fibered if the following two properties
hold:
\begin{itemize}
\item[(F1)] For each $y\in E$ and $\phi:b\to c=p(y)$ there exists a cartesian 
arrow $f:x\rTo y$ over $\phi$.
\item[(F2)] Composition of cartesian arrows is cartesian.
\end{itemize}

The choice of cartesian lifting for each arrow in $B$ is called {\em cleavage}.
Once such a choice is made, any arrow $\phi:b\to c$ in $B$ gives rise to an 
inverse image functor $\phi^*:E_c\rTo E_b$. This allows one to describe 
the fibered category $E$ as a pseudofunctor $B^\op\rTo\Cat$. If we are lucky,
the cartesian liftings chosen are compatible with the composition, and then
we get a strict functor from $B^\op$ to $\Cat$.

The inverse limit $\LIM(p)\in\Cat$ is defined as the category of cartesian 
sections
$s:B\rTo E$, that is the sections $s, \ p\circ s=\id_B$, carrying each arrow
of $B$ to a cartesian arrow. 

In applications from deformation theory we 
need to be able to make simultaneous weak groupoid completion of the 
fibers $E_b$. We get a collection of simplicial groupoids $\wt{E_b},\ b\in
B$. They are assembled into a quasifibered simplicial category over $B$,
see below.


\subsubsection{}
In what follows $[n]$ is the category defined by the ordered set
$\{0,\ldots,n\}$.
Let $B$ be a category. We denote by $\Delta^\op/B$ the following
category. Its objects are the functors
$\alpha:[n]^\op\rTo B$ in $B$.
\footnote{Of course, the categories $[n]$ and $[n]^\op$ are isomorphic.
We prefer the opposite category because of a notational convenience.}
A morphism
from $\beta:[m]^\op\to B$ to $\alpha:[n]^\op\to B$ is given by a functor 
$u:[m]\to[n]$ satisfying the condition $\beta=\alpha\circ u^\op$.

\begin{dfn}
A quasifibered category over $B$ is a (strict) functor 
$$ F: (\Delta^\op/B)^\op \rTo\Cat$$
satisfying the following property:
\begin{itemize}
\item[(QF)] For any $u:[m]\to[n]$ such that $u(0)=0$, $F(u)$ is an 
equivalence.  
\end{itemize}
\end{dfn}

\begin{rem*}A quasifibered category over $B$ should be considered as
a weak functor from $B^\op$ to the category of categories. The condition
(QF) which is formally of no use in this note, ensures that a map
$\bE\rTo\bE'$ of quasifibered categories over $B$ is an equivalence if
for each $b\in B$ the map of the fibers $\bE(b)\rTo\bE'(b)$ is an equivalence.
\end{rem*}

Let us show now how a fibered category $p:E\rTo B$ gives rise to a 
quasifibered category in our sense.
Let $\alpha:[n]^\op\rTo B$. We define 
$$\bE(\alpha)=\LIM(\pr_1:[n]^\op\times_BE\rTo\relax[n]^\op).$$
Thus, $\bE(\alpha)$ is the category of ``cartesian liftings'' of $\alpha$.

The property (QF) is fulfilled since the forgetful functor
$$ \LIM(\pr_1:[n]^\op\times_BE\to[n]^\op)\rTo E_{\alpha(0)}$$
is an equivalence: it is surjective on objects and bijective on morphisms.

Note once more that our definition avoids a choice of a cleavage.

\begin{lem}Let $p:E\rTo B$ be a fibered category. Then there is a natural
equivalence $\LIM(p)\rTo\holim(\bE)$.
\end{lem}
\begin{proof}
Comparing the definitions one immediately gets a canonical isomorphism
$\LIM(p)=\lim(\bE)$. One has a canonical map $\lim(\bE)\rTo\holim(\bE)$.
We prove it is an equivalence by checking that the functor 
$\bE:(\Delta^\op/B)^\op\rTo\Cat$
is fibrant in the Reedy model structure on the category of such functors,
see~\cite{hir}, Theorem 19.9.1(2), together with Propsition 15.10.4
and Definition 15.3.3.
This amounts to checking that for any $\alpha\in\Delta^\op/B$ the map
\begin{equation}\label{eq:matching}
\bE_\alpha\rTo M_\alpha\bE
\end{equation}
is a weak fibration in $\sCat$ where $M_\alpha$ is the matching object of 
$\bE$, see~\cite{hir}, Definition 15.2.3(2). If $\alpha$ is a $0$-simplex,
this follows from the fibrantness of the discrete categories. If $\alpha$
is an $n$-simplex with $n\geq 2$ then the map~(\ref{eq:matching}) is
a bijection. Finally, if $\alpha:b\to c$ is a $1$-simplex, the 
map~(\ref{eq:matching}) becomes
$$ \bE(\alpha)\rTo \bE(b)\times\bE(c).$$
Such map satisfies the conditions of Proposition~\ref{prop-weak-fibration} 
characterizing the weak fibrations.
\end{proof}

The notion of quasifibered category tautologically generalizes to any
category with a notion of weak equivalence instead of $\Cat$. For instance,
one has the notion of quasifibered simplicial category based on the weak
equivalence of simplicial categories.

The following claims are immediate.

\begin{prp}
\label{prp:old-nerve}
Let $\bE:(\Delta^\op/B)^\op\rTo\sCat$ be a quasifibered category. 
Then the 
composition of $\bE$ with the (standard) nerve functor $\cN$ 
$$ (\Delta^\op/B)^\op\rTo^\bE\sCat\rTo^{\cN}\simpl$$
is a quasifibered simplicial set.
\end{prp}
In fact, the standard simplicial nerve functor $\cN$ 
carries weak equivalences of simplicial categories to homotopy equivalences.
\qed

\begin{prp}
\label{prp:new-nerve}
Let $\bE:(\Delta^\op/B)^\op\rTo\sCat$ be a quasifibered category, 
so that for each $\sigma\in\Delta^\op/B$ the simplicial category $\bE_\sigma$
is fibrant.  Then the 
composition of $\bE$ with the nerve functor $\fN$ 
$$ (\Delta^\op/B)^\op\rTo^\bE\sCat\rTo^{\fN}\simpl$$
is a quasifibered simplicial set.
\end{prp}
In fact, the simplicial nerve functor $\fN$ carries weak equivalences
of fibrant simplicial categories to homotopy equivalences.
\qed

\begin{prp}
\label{prp:WQF}
Let $\bE:(\Delta^\op/B)^\op\rTo\sCat$ be a quasifibered simplicial category.
Then the composition of $\bE$ with the weak groupoid completion functor
$W\mapsto\wh{W}=L^H(W,W)$ (hammock localization with respect to all arrows) 
$$ (\Delta^\op/B)^\op\rTo^\bE\sCat\rTo^{L^H}\sCat$$
is a quasifibered simplicial category. In particular, fiberwise hammock 
localization of a fibered category is a quasifibered simplicial category.
\end{prp}
\begin{proof}
We have to show that a weak equivalence of simplicial categories
$f:\cC\to\cD$ induces a weal equivalence $\wh{f}$ of their weak groupoid
completions.

First of all, a weak equivalence
$f:\cC\rTo\cD$ of simplicial categories induces a homotopy equivalence
$\cN(f):\cN(\cC)\to\cN(\cD)$ of the classical nerves. Furthermore,
simplicial localization preserves the classical nerve, so one deduces
that $\cN(\wh{f})$ is an equivalence. Finally, since $\wh{f}$ is a morphism
of simplicial groupoids, \cite{dha}, 6.3.3, gives the required claim. 
\end{proof}


\begin{prp}\label{N-holim}
Let $\bE:(\Delta^\op/B)^\op\rTo\sCat$ be a quasifibered simplicial category
such that $\bE(\sigma)$ are fibrant groupoids
for all $\sigma\in\Delta^\op/B$. Then there is a natural 
homotopy equivalence
$$\fN(\holim(\bE))\sim \holim(\fN(\bE)).$$ 
\end{prp}
\begin{proof}
Since $\fN$ is right adjoint to $\CAT$, it preserves inverse limits. 
Since $\fN$ preserves strong equivalences and weak fibrations,
it carries simplicial framings to simplicial framings. 
\end{proof}

\section{Details}
\label{sec:details}

\subsection{The functor $S\mapsto\CAT(S)$}

The morphisms $\CAT^n(a,b)$ form simplicial ``standard cubes'' 
$(\Delta^1)^{b-a-1}$.

We present below a few useful notations.

\subsubsection{Standard cubes.}
Let $I=\Delta^1$ be the standard simplicial
segment. Its $k$-simplices are monotone functions 
$f:\{0,\ldots,k\}\rTo\{0,1\}$. We will denote $0$ and $1$ the corresponding
constant functions considered as $k$-simplices. These are $k$-simplices
obtained by degeneration from the zero simplices $0$ and $1$.

Let $\square=I^X$ where $X$ is a finite set. $k$-simplices of $\square$
are the collections $f=(f_x,\ x\in X)$ with $f_x\in I_k$.
The boundary $\p\square$ of $\square$
is defined as the union of $2|X|$ squares of dimension $|X|-1$
defined by the equations $f_x=0$ or $f_x=1$. Similarly to the
``simplicial horn'' $\Lambda^n_i$ obtained by erasing $i$-th face
from the boundary of the $n$-simplex, we define ``cubic horns''
$\Pi^n_{x,0}$ and $\Pi^n_{x,1}$ obtained by erasing from $\p\square$
the face given by the equation $f_x=0$ (resp., $f_x=1$).

\subsubsection{$\CAT^n$ are cofibrant.}
\label{cat-n-are-cofibrant}
Let $0\leq a< b\leq n$. A $k$-arrow
$f\in I^{(a,b)}$ is given by a collection $f=(f_{a+1},\ldots,f_{b-1})$
of its projections $[n]\rTo\left[1\right]$. 
There are $n+1$ such maps, according
to the number of preimages of $1$. We will write $f_i=1$ or $f_i=0$ 
for the respective constant maps.
Let 
$$\Ind(a,b)_m= \{f\in\shom(a,b)_m|f_x\ne 0\ \forall x\in(a,b)
\text{ and }f\not\in\cup_is_i(\shom(a,b)_{m-1}\}$$
be the collection of indecomposable nondegenerate $m$-arrows from $a$
to $b$. Put
$$ \Ind^n=\coprod_{m,0\leq a<b\leq n}\Ind(a,b)_m.$$

It is clear that $\Ind^n$ is the set of free generators of $\CAT^n$.
This proves that $\CAT^n$ is cofibrant.

\begin{prp}\label{inner-horns}
For $i\ne 0,n$ the map $\CAT(\Lambda^n_i)\rTo\CAT^n$
is a strong acyclic cofibration.
\end{prp}
\begin{proof}
An indecomposable map $f\in\CAT^n(a,b)_m$ belongs to the image
of 
$$\p^j: \CAT^{n-1} \rTo \CAT^n$$
 iff $j\not\in [a,b]$
or $j\in (a,b)$ and $f_j=1$. This implies that for $(a,b)\ne (0,n)$
$$\CAT(\Lambda^n_i)(a,b)=\CAT^n(a,b).$$
Let us identify the subset
$$\CAT(\Lambda^n_i)(0,n)\subseteq\CAT^n(0,n).$$
A simplex $f\in I^{(0,n)}_k,\ f=(f_1,\ldots,f_{n-1}),\ f_i\in I_k,$ 
belongs to the image iff there exists $j$ such that $f_j=0$ or there
exists $j\ne i$ with $f_j=1$. Thus, the image of 
$\CAT(\Lambda^n_i)(0,n)$ in $I^{(0,n)}$
is $\Pi^{(0,n)}_{i,1}$, the union of $2n-3$ faces (all but one) of 
$n-1$-dimensional cube. It is clearly contractible.
\end{proof}

\subsubsection{}
\label{lambda0n}
Let us explicitly describe the
map $\CAT(\Lambda^n_0)\rTo\CAT(\Delta^n)$.

Similarly to the above, the map
$$\CAT(\Lambda^n_0)(a,b)\rTo\CAT^n(a,b)$$
is a bijection for if $b<n$ or $a>1$. Thus, the non-trivial cases are
$(a,b)=(0,n)\text{ and }(1,n)$.

Let $(a,b)=(1,n)$. An element $f=(f_2,\ldots,f_{n-1})$
belongs to the image of $\CAT(\Lambda^n_0)(1,n)$ iff 
\begin{itemize}
\item either there exists $i\in(1,n)$ with $f_i=0$.
\item or there exists $i\in(1,n)$ with $f_i=1$.
\end{itemize}
This means that the image of $\CAT(\Lambda^n_0)(1,n)$ in $I^{(1,n)}$
is the boundary $\p I^{(1,n)}$.

Let $(a,b)=(0,n)$. An element $f=(f_1,\ldots,f_{n-1})$
belongs to the image of $\CAT(\Lambda^n_0)(0,n)$ if
\begin{itemize}
\item either there exists $i>1$ with $f_i=0$.
\item or there exists $i\geq 1$ with $f_i=1$.
\end{itemize}
Thus, the image of $\CAT(\Lambda^n_0)(0,n)$ in $I^{(0,n)}$
is the horn $\Pi^{(0,n)}_{1,0}$.

\subsubsection{Examples}
\label{cat-examples}
Let $S^n$ be the simplicial sphere defined as $S^n=\Delta^n/\p\Delta^n$.
The simplicial category $\CAT(S^n)$ has one object $*$. Therefore, $\CAT(S^n)$
is described by the simplicial monoid of the endomorphisms of $*$.
By definition, this is the free simplicial monoid generated by the set
$\{[\alpha]|\alpha\in\Mor\ \CAT^n\}$
 of
morphisms of $\CAT^n$, modulo the following relations:
\begin{itemize}
\item $[\alpha\beta]=[\alpha][\beta]$ if $\alpha$ and $\beta$ are composable.
\item $[\alpha]=1$ if $\alpha$ belongs to the image of $\p^j:\CAT^{n-1}
\rTo\CAT^n$ for some $j$.
\end{itemize}
Recall that in $\CAT^n$ one has $\Hom(a,b)=I^{(a,b)}$. A morphism
$f:a\to b$ belongs to the image of $\p^j$ iff $j\not\in (a,b)$ or
of $j\in (a,b)$ and $f_j=1$. Thus, the arrows $f:a\to b$ with $(a,b)\ne (0,n)$
disappear, as well as the arrows $f=(f_1,\ldots,f_{n-1})$ having some $f_i=1$
or $0$ (the latter means that $f$ is a nontrivial composition).

This proves that the simplicial monoid of endomorphisms describing $\CAT(S^n)$,
is free, generated by the simplicial set $S^1\wedge\ldots\wedge S^1$ 
($n-1$ times), where $S^1$ is considered as a pointed simplicial set
and $\wedge$ denoted the standard smash-product. 
The latter is, of course, homotopically $n-1$-sphere.

Let us now calculate $\CAT(D^{n+1})$ where 
$D^{n+1}=\Delta^{n+1}/\cup_{i>0}\p^i(\Delta^n).$

This simplicial category has as well one object.
Similarly to the above, the simplicial monoid describing 
$\CAT(D^{n+1})$ is the free
monoid generated by $[\alpha]$ with $\alpha\in\Mor\ \CAT(D^{n+1})$,
modulo relations
\begin{itemize}
\item $[\alpha]=1$ if $\alpha=\beta\gamma$ and $\beta:i\to n+1$ with $i>1$.
\item $[\alpha]=1$ if $\alpha$ belongs to the image of $\p^j:\CAT^{n-1}
\rTo\CAT^n$ for some $j>0$.
\end{itemize}

This implies that the simplicial monoid of endomorphisms describing 
$\CAT(D^{n+1})$,
is free, generated by the simplicial set $I^1\wedge\ldots\wedge S^1$ 
($n-1$ circle factors), where $I$ is considered as a pointed simplicial set 
with the marked point $1\in I$. This is homotopically $n$-disc.

\subsection{Proof of~\ref{prp:fibrations}}
\label{proof:fibrations}
Let $f:X\rTo Y$ be a fibration of fibrant groupoids. We have to check $f$
satisfies RLP with respect to the maps $\CAT(\Lambda^n_i)\rTo \CAT(\Delta^n)$.

The case $i\ne 0,n$ follows from~\ref{inner-horns}.

\subsubsection{}
We will start with the first non-trivial case to make the general
case presented in~\ref{proof:fib-gen} more readable. 

Let $n=2,\ i=0$. Thus, we have the 
following data: 
\begin{itemize}
\item The objects $x_0,x_1,x_2$ of $X$ with $0$-morphisms $u_{01}:x_0\to x_1$
and $u_{02}:x_0\to x_2$.
\item The images of the above $y_i=f(x_i), v_{0i}=f(u_{0i})$.
\item A $0$-morphism $v_{12}:y_1\to y_2$ and $\beta\in\shom_1(y_0,y_2)$
such that $$\beta:v_{12}\cdot v_{01}\rTo v_{02}.$$
\end{itemize}

We have to lift $v_{12}$ and $\beta$.

{\em Step 1.} First of all, since $X$ is a fibrant groupoid, there exists
$u:x_1\to x_2$ and a one-morphism $\theta:u\cdot u_{01}\rTo u_{02}$.

{\em Step 2.} Since $\shom(y_0,y_2)$ is fibrant, the pair of $1$-simplices
$\beta$ and $f(\theta)$ can be completed to a $2$-simplex $\Xi$; its
remaining edge will be denoted $\eta:f(u)\cdot v_{01}\rTo v_{12}\cdot v_{01}$.

$$
\begin{diagram}
& & v_{12}\cdot v_{01} & & \\
& \ruTo^\eta & & \rdTo^{\beta}&\\
f(u)\cdot v_{01} & &\rTo^{f(\theta)} & & v_{02} 
\end{diagram}
$$

{\em Step 3.} We wish to ``correct'' now $\Xi$ so that $\eta$ will be of form
$\bar{\eta}\cdot v_{01}$ for some $\bar{\eta}:f(u)\to v_{12}$.

This is possible for the following reason. The right composition
$$ \shom(y_1,y_2)\rTo^{\cdot v_{01}}\shom(y_0,y_2)$$
is a weak equivalence of Kan simplicial sets. We have a pair of
points $f(u)$ and $v_{12}$ in the left-hand side, and we have a path
$\eta$ connecting their images in the right-hand side. This implies there
is a path $\bar{\eta}$ on the left-hand side whose image is homotopic to 
$\eta$. The last thing to do is to use the mentioned above homotopy
to correct $\Xi$.

Thus, we assume in the sequel that $\eta=\bar{\eta}\cdot v_{01}$.

{\em Step 4.} Since $\shom(x_1,x_2)\rTo\shom(y_1,y_2)$ is a fibration,
the $1$-simplex $\bar{\eta}:f(u)\rTo v_{12}$ lifts to a $1$-simplex
$\tilde{\eta}:u\rTo u_{12}$ (thus defining $u_{12}$).

{\em Step 5.} Finally, since $\shom(x_0,x_2)\rTo\shom(y_0,y_2)$ is a 
fibration,
the $2$-simplex $\Xi$ downstairs lifts to a $2$-simplex
$$
\begin{diagram}
& & u_{12}\cdot u_{01} & & \\
& \ruTo & & \rdTo^{\alpha}&\\
u\cdot u_{01} & &\rTo & & u_{02} 
\end{diagram}
$$
In particular, this assures that $f(\alpha)=\beta$ and this completes the
proof in the case $n=2,\ i=0$.

\subsubsection{}
\label{proof:fib-gen}

Let now $n>2$. The steps below follow the steps 1 to 5  of the special case 
$n=2$.

It is worthwhile to use the diagram~(\ref{kabir}) below to follow
the construction of the simplicial sets and of the arrows. Smaller
``local'' diagrams are presented in the text.

{\em Step 1.} Let us check that the map $\CAT(\Lambda^n_0)\rTo X$ extends
to a map $\CAT^n\rTo X$. This amounts to say that the following commutative
diagram
$$
\begin{diagram}
\p I^{(1,n)} & \rTo & \shom(x_1,x_n) \\
\dTo & & \dTo_{\cdot u_{01}} \\
\Pi^{(0,n)}_{1,0} & \rTo & \shom(x_0,x_n)
\end{diagram},
$$
the right arrow being the composition with the arrow 
$u_{01}:x_0\rTo x_1$, extends to the 
following diagram
$$
\begin{diagram}
\p I^{(1,n)} & \rTo^{i_1}& I^{(1,n)} & \rTo^{\rho_1} & \shom(x_1,x_n) \\
\dTo^{j_1} & &\dTo^{j_2} & & \dTo_{\cdot u_{01}} \\
\Pi^{(0,n)}_{1,0} &\rTo^{i_2} & I^{(0,n)}& \rTo^{\rho_2} & \shom(x_0,x_n)
\end{diagram}.
$$
The left-hand commutative square in the diagram is equivalent
to presentation of the boundary of a $(n-1)$-dimensional disc as a union 
of two hemispheres along a common $(n-3)$-sphere.

Since the composition with $u_{01}$ is a homotopy equivalence of Kan 
simplicial sets, the claim easily follows.


{\em Step 2.} The map $\CAT^n\rTo X$ composed with $f:X\rTo Y$, does not give,
of course, the given map  $\CAT^n\rTo Y$. However, the two maps
coincide on $\CAT(\Lambda^n_0)$. 
This gives the following commutative diagram in $\simpl$
$$
\begin{diagram}
\Pi^{(0,n)}_{1,0} &\rTo^{i_2}& I^{(0,n)}& \rTo & \shom(x_0,x_n) \\
\dTo^{k_2} & & & & \dTo^{f_{0n}} \\
I^{(0,n)} & \rTo&  &  & \shom(y_0,y_n)
\end{diagram}.
$$
 Since $\shom(y_0,y_n)$ is fibrant, 
this diagram can be completed to the following
\begin{equation}
\begin{diagram}\label{diag:Xi}
\Pi^{(0,n)}_{1,0} &\rTo^{i_2} & I^{(0,n)}& \rTo & \shom(x_0,x_n) \\
\dTo^{k_2} & &\dTo^{k_3} & & \dTo^{f_{0n}} \\
I^{(0,n)} & \rTo^{i_3}& \Xi &\rTo^\xi  & \shom(y_0,y_n)
\end{diagram},
\end{equation}
where $\Xi$ is defined by the cocartesian diagram
$$
\begin{diagram}
I\times\Pi^{(0,n)}_{1,0} &\rTo & I\times I^{(0,n)} \\
\dTo^{\pr_2} &  & \dTo \\
\Pi^{(0,n)}_{1,0} & \rTo &\Xi
\end{diagram}.
$$

In fact, the left-hand side commutative square of~(\ref{diag:Xi}) 
induces an acyclic cofibration
$$ I^{(0,n)}\coprod^{\Pi^{(0,n)}_{1,0}}I^{(0,n)}\rTo \Xi.$$
Note that the boundary of $\Xi$ (i.e. the image of the boundary of
$I\times I^{(0,n)}$ under the natural projection) consists
of two copies of $I^{(0,n)}$ and of the image $\mathrm{H}$ of 
$I\times I^{(1,n)}$
\footnote{$\mathrm{H}$ is the capital Eta}
which is defined by the following cocartesian diagram
 $$
\begin{diagram}
I\times\p I^{(1,n)} &\rTo & I\times I^{(1,n)} \\
\dTo^{\pr_2} &  & \dTo \\
\p I^{(1,n)} & \rTo &\mathrm{H}
\end{diagram}.
$$

{\em Step 3.} We wish to correct $\Xi$ so that the restriction $\eta$ of
$\xi:\Xi\rTo\shom(y_0,y_n)$ (see~(\ref{diag:Xi}) and (\ref{kabir})) to 
$\mathrm{H}$ 
is the composition
of a map $\bar{\eta}:\mathrm{H}\rTo\shom(y_1,y_n)$ and of the map 
$v_{01}:y_0\rTo y_1$. This is possible for the same reason as in the case 
$n=2$.

Thus, from now on we assume that $\eta=\bar{\eta}\cdot v_{01}$.
 
{\em Step 4.} Since $\shom(x_1,x_n)\rTo\shom(y_1,y_n)$ is a fibration,
and since the square consisting of the arrows $\rho_1,\ \rho_2,\ j_2$ and
$\cdot u_{01}$ is commutative (see~(\ref{kabir}) below), there exists
a dashed arrow $\delta:\mathrm{H}\rTo\shom(x_1,x_n)$ making the diagram 
commutative.

{\em Step 5.} One has now two maps, $\rho_2:I^{(0,n)}\rTo\shom(x_0,x_n)$
and \newline
$\delta:\mathrm{H}\rTo\shom(x_1,x_n)$ whose compositions with $j_2$
and with $k_4$ correspondingly, coincide. This defines a map
$$\mathrm{H}\coprod^{I^{(1,n)}}I^{(0,n)}\rTo\shom(x_0,x_n).$$
It extends to a map $\Delta:\Xi\rTo\shom(x_0,x_n)$ since the map
$$\mathrm{H}\coprod^{I^{(1,n)}}I^{(0,n)}\rTo\Xi$$
is an acyclic cofibration.
  
Look now at the following diagram we have constructed
\footnote{
we have not found a place to draw the arrow
$\Delta:\Xi\rTo\shom(x_0,x_n)$, as well as 
to name the dashed arrow by $\delta$.
}

\begin{equation}\label{kabir}
\begin{diagram}[small]
\p I^{(1,n)}&     &\rTo^{i_1}        &       &I^{(1,n)}&          
             &\rTo^{\rho_1}&&    &\shom(x_1,x_n)  &         &   \\
            &\rdTo^{j_1}&            &       &\vLine^{k_4}&\rdTo^{j_2}&        
            &&\ruDashto(5,4)&\vLine^{f_{1n}}&\rdTo^{\cdot u_{01}} &   \\
\dTo^{k_1}     &     &\Pi^{(0,n)}_{1,0} &\rTo^{i_2}&\HonV &          
&I^{(0,n)}&&\rTo^{\rho_2}&\HonV         &         &\shom(x_0,x_n)   \\
            &     &\dTo^{k_2}        &       &\dTo^{}&         &\dTo_{k_3}  
&&       &\dTo^{}       &                  &\dTo^{f_{0n}}  \\
I^{(1,n)}   &\rLine^{i_4} &\VonH           &\rTo^{}&\mathrm{H}  &\rLine    
&\VonH    &&\rTo^{\bar{\eta}\hspace{1cm}}&\shom(y_1,y_n)& &  \\
            &\rdTo&                  &       &         &\rdTo^{}  &         
&&     &                &\rdTo^{\cdot v_{01}} &  \\
            &     &I^{(0,n)}         &\rTo^{i_3}&      &          &\Xi      
&&\rTo^\xi&             &         &\shom(y_0,y_n)  \\
\end{diagram}
\end{equation}

The composition $\delta\cdot i_4: I^{(1,n)}\rTo \shom(x_1,x_n)$, together 
with the composition
$\Delta\cdot i_3: I^{(0,n)}\rTo \shom(x_0,x_n)$, define a functor 
$\CAT^n\rTo X$ extending
the given functor $\CAT(\Lambda^n_0)\rTo X$; this is what we had to do.

\subsection{Examples}We present here two counterexamples showing that 
fibrations are not always preserved by the nerve functors.
\label{ss:examples}
\subsubsection{The classical nerve $\cN$.}
\label{sss:cN-ex}
We present below an example of a fibration of fibrant groupoids $\cC\rTo\cD$
for which the map $\cN(\cC)\rTo\cN(\cD)$ is not Kan fibration.
All simplicial categories in the example will have three objects $0,1,2$,
with the functors identical on objects.

Choose $\cD$ to be the contractible (discrete) groupoid with three objects. 
The standard triangle $\SC^2$ has three objects $0,1,2$ with morphisms
generated by the arrows $v_{01}\in\shom_2(0,1)$ and $v_{12}\in\shom_2(1,2)$.
The $0$-horn $\SC(\Lambda^2_0)$ has the same objects, and the morphisms
freely generated by $u_{01}\in\shom_1(0,1)$ and $u_{02}\in\shom_1(0,2)$.
The embedding of $\SC(\Lambda^2_0)$ to $\SC^2$ carries $u_{01}$ to $d_2v_{01}$
and $u_{02}$ to $d_1v_{12}\circ d_1v_{01}$. The composition
$$\SC(\Lambda^2_0)\rTo\SC^2\rTo \cD$$
can be decomposed as
$$\SC(\Lambda^2_0)\rTo\cC\rTo \cD,$$
the first map being a cofibration, and the second an acyclic fibration.
Thus, $\cC$ is automatically a fibrant groupoid. The image of $d_1u_{02}$ 
in $\cC$ can not be presented as a composition of
$d_1u_{01}$ with another arrow. Therefore, the map of $\SC^2$  to $\cD$
can not be lifted to $\cC$.

\subsubsection{The homotopy coherent nerve $\fN$.}
Another example we wish to present is that of a strong fibration of fibrant
simplicial categories $\cC\rTo\cD$, so that the map $\fN(\cC)\rTo\fN(\cD)$
is not a Kan fibration. Here $\cC$ and $\cD$ are discrete categories
with the same objects $0,1,2$; $\cD$ has three arrows $v_{01},\ v_{12}$
and $v_{02}=v_{12}\circ v_{01}$. The category $\cC$ has the same objects,
with the arrows $u_{ij}$ so that $u_{02}=u_{12}\circ u_{01}$, and 
$w_{02}:0\to 2$. The map $\cC\rTo \cD$ identical on objects is uniquely defined
and it is a fibration.

Note that $\CAT(\Lambda^2_0)$ has three objects $0,1,2$ and two arrows from
$0$ to $1$ and to $2$. The simplicial category $\CAT^2$ has one more $0$-arrow
from $1$ to $2$ and one more $1$-arrow. Define the functor $\CAT^2\rTo\cD$
in the only possible way; if we define a functor  $\CAT(\Lambda^2_0)\rTo\cC$
so that the arrow $w_{02}$ is in the image, we cannot complete the 
corresponding commutative diagram.

\section{Errata et Corrigenda for \cite{dgc,dha,dsa}}
\label{sec:corr}

In this section we indicate what should be fixed in the papers 
\cite{dgc,dha,dsa}. 

\subsection{General corrections}

In the context of~\cite{dgc,dha,dsa} deformation functors
appear as functors  $\sDef$ with values in simplicial categories.
Another avatar of deformation functor, $\Def$, is defined as the composition of
$\sDef$ with the nerve functor $\cN$.
 
In order to satisfy the announced properties, one has to replace the
nerve functor $\cN$ with $\fN$ in the definition of $\Def$. Thus,
we replace the old definition with the following.
\begin{dfn}
$$\Def=\fN\circ\sDef.$$
\end{dfn}

Both definitions coincide up to homotopy since $\sDef$ has values in 
fibrant groupoids as follows from the lemma below.

\begin{lem}\label{wstar-is-fibrant}
Let $\cC$ be a  model category endowed with a functorial simplicial
framing.
Then the simplicial category 
$W_*(\cC)$ whose objects are the fibrant cofibrant objects of $\cC$,
and degree $n$ maps are weak equivalences $x\rTo y^{\Delta^n}$, is a fibrant
groupoid.
\end{lem}
\begin{proof}Obviously $\pi_0(W_*(\cC))$ is a groupoid. If $x,y\in W_*(\cC)$,
the simplicial set $\shom(x,y)$ is fibrant by the properties of the 
simplicial framing.
\end{proof}


\subsubsection{}
The important descent property for the functor $\Def$ results from the 
descent property of the functor $\sDef$ and the following result 
(see Proposition~\ref{N-holim}) which replaces~\cite{dgc}, 
A5.2.

\begin{prp*}
Let $\bE:(\Delta^\op/B)^\op\rTo\sCat$ be a quasifibered simplicial category
such that $\bE(\sigma)$ are fibrant groupoids
for all $\sigma\in\Delta^\op/B$. Then there is a natural 
homotopy equivalence
$$\fN(\holim(\bE))\sim \holim(\fN(\bE)).$$ 
\end{prp*}

Descent property for the functor $\sDef$ formulated in \cite{dsa}, 3.2.3,
should be now reformulated in the context of quasifibered categories
as follows.

Let $V_\bullet$ be a hypercovering of $X$ and let $\cO$ be an operad over $X$.
Recall that for a commutative algebra $R$ \ $\cW_\sim^\fl(R,X)$ denotes the 
category of sheaves of $\cO$-algebras on $X$ over $R$, flat over $R$,
with arrows being quasiisomorphisms of algebras.

The hypercovering $V_\bullet$ of $X$ defines a fibered category
$\cW_\sim^\fl(R,V_\bullet)$ over $\Delta^\op$, with the fiber
$\cW_\sim^\fl(R,V_n)$ over $n\in\Delta^\op$. By \ref{prp:WQF} the full 
Dwyer-Kan localization yields a quasifibered simplicial category 
$$\bbW:(\Delta^\op/\Delta^\op)^\op\rTo \sCat,\quad 
\bbW(n)=\wh{\cW}_\sim^\fl(R,V_n).$$

The canonical map
\begin{equation}
\label{eq:descent}
 \wh{\cW}_\sim^\fl(R,X)\rTo\holim\ \bbW
\end{equation}
is obtained from the comparison of $V_\bullet$ to the trivial hypercovering
of $X$.

The proof that (\ref{eq:descent}) yields the descent property for the 
deformation functor, given in Section~3.2.3 of~\cite{dsa}, is left untouched
(the nerve functor $\cN$ being replaced with $\fN$).

\subsection{Errata for \cite{dgc}}

Proposition A3.2 of ~\cite{dgc} claims that the standard nerve functor
$\cN$ preserves the weak equivalences, the cofibrations and the fibrations.
The first two claims are corect; the claim about the fibrations is wrong,
see \ref{sss:cN-ex}. The results of~\ref{ss:newnerveprop} describing the
properties of the new nerve $\fN$ replace the wrong statement.

Lemma A3.3 (used in A4.3) is wrong. The present Section 3
describing a choice of simplicial framing 
together with the notion of quasifibered category, replaces subsection A4.3.

\subsection{Errata for~\cite{dha}}

Proposition 6.3.3 remains correct (and useful) as it is formulated (for the
old nerve functor $\cN$). Since for fibrant groupoids the nerves $\cN$ and
$\fN$ are canonically equivalent, it remains true also with the new
simplicial nerve. Thus, one has

\begin{prp}
\label{converse-we}
Let $f:\cC\rTo\cD$ be a map of fibrant groupoids. If $\fN(f)$
is a weak equivalence, $f$ is as well a weak (and strong) equivalence.
\end{prp}

\subsection{Errata for~\cite{dsa}}

Sections A1.9-A1.11 of~\cite{dsa} (simplical structure, pseudofunctors and 
homotopy limits) are wrong or obsolete. See Section~\ref{sec:qf} instead.

A1.12 (which is ~\cite{dgc}, A5.2) is replaced with 
Proposition~\ref{N-holim}.

In formula A1.13 \ $\cD^{\Delta^1}$ should be meant in the sense of
the chosen simplicial framing, see~\ref{ss:holim}.

Corollary A2.3  is correct with the standard version $\cN$ of the nerve.
For the new version of nerve $\fN$ one has Proposition~\ref{converse-we}
instead.

Since there is no canonical simplicial structure on $\sCat$, claim A.4 of
\cite{dsa} has no meaning. In applications Proposition~\ref{prp:WQF} 
is used to replace it.

Corollary A2.5 is correct as it is, with the standard vesion of the nerve.

\end{document}